\documentclass[12pt]{article}
\usepackage{amsfonts}
\usepackage{amssymb}
\usepackage{amsmath}
\usepackage{eufrak}
\usepackage[dvips]{graphicx}
\usepackage{latexsym}
\newcommand{\ind}{\makebox[1em]{\raisebox{-.5ex}[0ex][0ex]{\makebox[0em]%
{$\smile$}}\raisebox{.4ex}[0ex][0ex]{\makebox[-.02em]{$|$}}}}

\newcommand{\dep}{\makebox[1em]{\raisebox{.3ex}[0ex][0ex]%
{$\not$}\makebox[.7em]{\ind}}}

\newcommand{\nmdep}{\makebox[1em]{\raisebox{1.5ex}[0ex][0ex]{\makebox[0em]%
{$\scriptscriptstyle n\! m$}}\makebox[-1em]{$\dep$}}}
\newcommand{\nmind}{\makebox[1em]{\raisebox{1.5ex}[0ex][0ex]{\makebox[0em]%
{$\scriptscriptstyle n\! m$}}\makebox[-1em]{$\ind$}}}

\newcommand{\bi}{\begin{itemize}}
\newcommand{\ei}{\end{itemize}}

\newcommand\NM{\mathcal{N\!M}}
\newtheorem{theorem}{Theorem}[section]

\newtheorem{fact}[theorem]{Fact}

\newtheorem{corollary}[theorem]{Corollary}
\newtheorem{proposition}[theorem]{Proposition}
\newtheorem{definition}[theorem]{Definition}

\newtheorem{question}[theorem]{Question}

\newtheorem{problem}[theorem]{Problem}

\title{New examples of small Polish structures}
\author{Jan Dobrowolski}
\date{}

\begin{document}
\maketitle
\begin{abstract}We answer some questions from \cite{1} by giving suitable examples of small Polish structures. First, we present a class of small Polish group structures without generic elements. Next, we construct a first example of a small non-zero-dimensional Polish $G$-group.
\end{abstract}

\section{Introduction}
In \cite{1}, Krupi\'{n}ski defined and investigated Polish structures by methods motivated by model theory.
\begin{definition}
A Polish structure is a pair $(X,G)$, where $G$ is a Polish group acting faithfully on a set $X$ so that the stabilizers of all singletons are closed subgroups of $G$. We say that $(X,G)$ is small if for every $n<\omega$, there are only countably many orbits on $X^n$ under the action of $G$.
\end{definition}
A particularly interesting situation is when the underlaying set $X$ is a group itself.
Throughout this paper, we follow the terminology from \cite{1}.
\begin{definition}
Let $G$ be a Polish group. \\
(i) A Polish group structure is a Polish structure $(H,G)$ such that $H$ is a group and $G$ acts as a group of automorphisms of $H$. \\
(ii) A (topological) $G$-group is a Polish group structure $(H,G)$ such that $H$ is a topological group and the action of $G$ on $H$ is continuous. \\
(iii) A Polish [compact] $G$-group is a topological $G$-group $(H,G)$, where $H$ is a Polish [compact] group. 
\end{definition}
Let $(X,G)$ be a Polish structure. For any finite $C\subseteq X$, by $G_C$ we denote the pointwaise stabilizer of $C$ in $G$, and for a finite tuple $a$ of elements of $X$, by $o(a/C)$ by denote the orbit of $a$ under the action of $G_C$ (and we call it the orbit of $a$ over $C$).

A fundamental concept for \cite{1} is the relation of $nm$-independence in an arbitrary Polish structure.
\begin{definition}
 Let $a$ be a finite tuple and $A$, $B$ finite subsets of $X$. Let $\pi _A:G_A\to o(a/A)$ be defined by $\pi _A(g)=ga$. We say that $a$ is $nm$-independent from $B$ over $A$ (written $a\nmind_A B$) if $\pi _A^{-1}[o(a/AB)]$ is non-meager in $\pi _A^{-1}[o(a/A)]$. Otherwise, we say that $a$ is $nm$-dependent on $B$ over $A$ (written $a\nmdep_A B$).
\end{definition}

This is a generalization of $m$-independence, which was introduced by Newelski for profinite structures.
Under the assumption of smallness $nm$-independence has similar properties to those of forking independence in stable theories, and hence it allows to transfer some ideas and techniques from stability theory to small Polish structures(which are purely topological objects). The investigation of Polish structures has been undertaken in \cite{4} and \cite{5}. For example, in \cite{4}, some structural theorems about compact $G$-groups were proved, and in \cite{5}, dendrites were considered as Polish structures, and some properties introduced in \cite{1} were examined for them.

The class of Polish structures contains much more interesting examples from the classical mathematics than the class of profinite structures. For example, for any compact metric space $P$, if we consider the group $Homeo(P)$ of all homeomorphisms of $P$ equipped with the compact-open topology, then $(P,Homeo(P))$ is a Polish structure (examples of small Polish structures of this form were investigated in \cite{1} and \cite{5}).
However, in the class of small Polish group structures, it is more difficult to construct interesting examples. In the present paper, we answer some questions from \cite{1} by constructing suitable examples of small Polish group structures.

The following is \cite[Question 5.4]{1} (see Definition \ref{def:gen} for the notion of $nm$-generic orbit).
\begin{question}\label{qgeneric}
Let $(G,H)$ be a small Polish group structure. Does $H$ possess an $nm$-generic orbit?
\end{question}
Proposition 5.5 from \cite{1} gives us a positive answer to Question \ref{qgeneric} in the class of small Polish $G$-groups.
In Section 2, we construct a class of small Polish group structures for which the answer to Question \ref{qgeneric} is negative.

The following problem was formulated in \cite{1} (after Question 5.32):
\begin{problem}\label{dimension}
Find a non-zero-dimensional, small Polish $G$-group.
\end{problem}
In Section \ref{dim}, we construct a small Polish $G$-group $(H,G)$, such that $H$ is homeomorphic to the complete Erd\"{o}s space, which is known to be one-dimensional.

\section{Preliminaries} 
If $A$ is a finite subset of $X$ (where $(X,G)$ is a Polish structure), we define the algebraic closure of $A$ (written $Acl(A)$) as the set of all elements of $X$ with countable orbits over $A$. If $A$ is infinite, we define $Acl(A)=\bigcup\{Acl(A_0):A_0\subseteq A$ is finite$\}$.
By Theorems 2.5 and 2.10 from \cite{1}, we have:
\begin{theorem}\label{properties}
In any Polish structure $(X,G)$, $nm$-independence has the following properties:\\
(0) (Invariance) $a\nmind_A B\iff g(a)\nmind_{g[A]} g[B]$ whenever $g\in G$ and $a,A,B\subseteq X$ are finite.\\
(1) (Symmetry) $a\nmind_C b\iff b\nmind_C a$ for every finite $a,b,C\subseteq X$.\\
(2) (Transitivity) $a\nmind_B C$ and $a\nmind_A B$ iff $a\nmind_A C$ for every finite $A\subseteq B\subseteq C\subseteq X$ and $a\subseteq X$.\\
(3) For every finite $A\subseteq X$, $a\in Acl(A)$ iff for all finite $B\subseteq X$ we have $a\nmind_A B$.\\
If additionally $(X,G)$ is small, then we also have:\\
(4) (Existence of $nm$-independent extensions) For all finite $a\subseteq X$ and $A\subseteq B\subseteq X$ there is $b\in o(a/A)$ such that $b\nmind_A B$.
\end{theorem}
The notion of $nm$-independence leads to the definition of $\NM$-rank.
\begin{definition}
The $\NM$-rank is the unique function from the collection of orbits over finite sets to the ordinals together with $\infty$, satisfying\\  $\NM(a/A)\geq\alpha +1$ iff there is a finite set $B\supseteq A$ such that $a\nmdep_A B$ and $\NM(a/B)\geq\alpha$. \\ The $\NM$-rank of $X$ is defined as the supremum of $\NM (x/\emptyset)$, $x\in X$.\\
We say that a Polish structure $(X,G)$ is $nm$-stable, if $\NM(X)<\infty$.
\end{definition}
The following fact is a part of \cite[Proposition 2.3]{1}.
\begin{fact}\label{indep}
Let $(X,G)$ be a Polish structure, $a$ be a finite tuple and $A,B$ be finite subsets of $X$. Then, TFAE:\\
(1) $a\nmind_A B$ \\
(2) $G_{AB}G_{Aa}\subseteq_{nm}G_A$
\end{fact}
By \cite[2.14]{1}, under some assumptions, $nm$-dependence in a $G$-group $(G,H)$ can be expressed in terms of the topology on $H$:
\begin{theorem}\label{gdelta}
Let $(X,G)$ be a Polish structure such that $G$ acts continuously on a Hausdorff space $X$. Let $a,A,B\subseteq X$ be finite. Assume that $o(a/A)$ is non-meager in its relative topology. Then, $a\nmind_A B \iff o(a/AB)\subseteq _{nm}o(a/A)$.
\end{theorem}
Counterparts of various notions from model theory were studied by Krupi\'{n}ski in the context of Polish structures. One of them is the notion of a generic orbit:
\begin{definition}\label{def:gen}
Let $(H,G)$ be a Polish group structure.
We say that the orbit $o(a/A)$ is left $nm$-generic (or that $a$ is left $nm$-generic over $A$) if for all $b\in H$ with $a\nmind_A b$, one has that $b\cdot a\nmind A,b$. We say that it is right $nm$-generic if, for $b$ as above, we have $a\cdot b\nmind A,b$. An orbit is $nm$-generic if it is both right and left $nm$-generic. 
\end{definition}
It was noticed in \cite{1} that $nm$-generics have similar properties to generics in simple theories, e.g. being right $nm$-generic is equivalent to being left $nm$-generic. We recall Proposition 5.5 from \cite{1}, which gives us a positive answer to Question \ref{qgeneric} for the class of small $G$-groups $(H,G)$ in which $H$ in not meager in itself (this holds, for example, in all Polish $G$-groups).
\begin{fact}\label{generic}
Suppose $(H,G)$ is a small $G$-group. Assume $H$ is not meager in itself (e.g. $H$ is Polish or compact, or, more generally, Baire). Then, at least one $nm$-generic orbit in $H$ exists, and an orbit is $nm$-generic in $H$ iff it is non-meager in $H$.
\end{fact}

Consider any $p\geq 1$ and the Banach space $l_p$ (over $\mathbb{R}$). We extend the $p$-norm from $l_p$ to $\mathbb{R}^{\omega}$ by putting $||z||=\infty$ for every $z\in \mathbb{R}^{\omega}\backslash l_p$.
The complete Erd\"{o}s space is the intersection of $l_p$ with $(\mathbb{R}\backslash \mathbb{Q})^{\omega}$ (with the topology induced from $l_p$).

Let $E_0,E_1,\dots$ be a fixed sequence of subsets of $\mathbb{R}$ and put 
$$\varepsilon =l_p\cap\prod_{n<\omega}E_n.$$
The following is a part of \cite[Theorem 1]{2}:
\begin{theorem}\label{criterion}
Assume that $\varepsilon$ is not empty and that every $E_n$ is zero-dimensional. For each $k\in \omega\backslash\{0\}$ we let $\eta(k)\in\mathbb{R}^{\omega}$ be given by
$$\eta(k)_n=sup\{|a|:a\in E_n\cap [-1/k,1/k]\},$$
where $\sup\emptyset=0$. The following statements are equivalent: \\
(1) $||\eta(k)||=\infty$ for each $k\in\omega\backslash\{0\}$ \\
(2) $\dim \varepsilon>0$
\end{theorem}
Under the assumptions of the above theorem, also the following theorem was proved there (\cite[Theorem 3]{2}):
\begin{theorem}\label{erdos}
If every $E_i$ is closed in $\mathbb{R}$, then $\varepsilon$ is homeomorphic to the complete Erd\"{o}s space if and only if $\dim\varepsilon>0$ and every $E_n$ is zero dimensional.
\end{theorem}

\section{Small Polish group structures without generic elements}
In this section, we construct a class of small Polish group structures for which the answer to Question \ref{qgeneric} is negative.

Suppose $(X,G)$ is a Polish structure. Let $H$ be an arbitrary group. For any $x\in X$ we consider an isomorphic copy $H_x=\{h_x:h\in H\}$ of $H$. By $H(X)$ we will denote the group $\bigoplus_{x\in X} H_x$. Although $H(X)$ is not necessarily commutative, we will use its group action only for commuting elements and we denote it by $+$. For any $y\in H(X)$ there are $h_1,\dots,h_n\in H\backslash\{e\}$ and pairwaise distinct $x_1,\dots,x_n\in X$ such that $y=(h_1)_{x_1}+\dots +(h_n)_{x_n}$. Then, by $\tilde{y}$ we will denote the set $\{x_1,\dots,x_n\}$. We also put $\tilde{A}=\bigcup_{y\in A}\tilde{y}$ for any $A\subseteq H(X)$.

Group $G$ acts as automorphisms on $H(X)$ by $$g((h_1)_{x_1}+\dots +(h_n)_{x_n})=(h_1)_{gx_1}+\dots +(h_n)_{gx_n}.$$ 
It is easy to see that if $h_1,\dots ,h_k\in H$ are pairwaise distinct, and\\ $x_{1,1},\dots,x_{1,i_1},x_{2,1},\dots,x_{2,i_2},\dots,x_{k,1},\dots,x_{k,i_k}\in X$ are pairwaise distinct as well, then \begin{equation}\begin{split}  G_{(h_{1})_{x_{1,1}}+\dots +(h_{1})_{x_{1,i_1}}+\dots +(h_{k})_{x_{k,1}}+\dots +(h_{k})_{x_{k,i_k}}}= \\
 =\bigcup_{\sigma_1\in S_{i_1},\dots ,\sigma_k\in S_{i_k}}\bigcap_{l\in\{1,\dots,k\}}\bigcap_{j\in\{1,\dots,i_l\}}\{g\in G:gx_{l,j}=x_{l,\sigma_l(j)}\}. \label{1} \end{split}\end{equation} By (\ref{1}), we get that for every $a\in H(X)$, $G_{\tilde{a}}$ is a subgroup of finite index in $G_a$, and hence, for every finite $A\subseteq H(X)$, $G_{\tilde{A}}$ is a subgroup of finite index in $G_A$.

\begin{proposition}
If $(X,G)$ is a Polish structure, and $H$ is a group, then $(H(X),G)$ is a Polish group structure. If, additionally, $(X,G)$ is small and $H$ is countable, then $(H(X),G)$ is small.
\end{proposition}
{\em Proof.} For any $a\in H(X)$ we have that $G_{\tilde{a}}$ is closed in $G$ and has finite index in $G_{a}$, so, $G_a$ is also closed in $G$. Hence, $(H(X),G)$ is a Polish group structure.

Now, assume that $(X,G)$ is small, and $H$ is countable. Then, for every fixed $k<\omega$ and $i_{1},\dots, i_{n}<\omega$, the orbit of a tuple $((h_{1,1})_{x_{1,1}}+\dots +(h_{1,i_1})_{x_{1,i_1}},\dots ,(h_{k,1})_{x_{k,1}}+\dots +(h_{k,i_k})_{x_{k,i_k}})$ depends only on $h_{1,1},\dots,h_{1,i_1},\dots,h_{k,1},\dots,h_{k,i_k}$ and on the orbit of the tuple $(x_{1,1},\dots,x_{1,i_1},\dots,x_{k,1},\dots,x_{k,i_k})$ in $(X,G)$. So, there are only countably many $k$-orbits in $(H(X),G)$. \hfill $\square$
\begin{proposition} \label{new ind}
Let $(X,G)$ be a Polish structure, and $H$ a countable group. Then for any finite $A,B,C\subseteq H(X)$, we have:\\
(1) $A\nmind_C B \iff \tilde{A}\nmind_{\tilde{C}} \tilde{B}$.\\
(2) If $a$ is a finite tuple of elements of $H(X)$, and $b$ is a tuple of elements of $X$ enumerating $\tilde{a}$,
 then $\NM(a/A)=\NM(b/\tilde{A})$. In particular, $(H(X),G)$ is $nm$-stable iff $(X,G)$ is.
\end{proposition}
{\em Proof.} (1) By Fact \ref{indep}, it is enough to show that $G_{CB}G_{CA}\subseteq_{nm}G_C \iff G_{\tilde{C}\tilde{B}}G_{\tilde{C}\tilde{A}}\subseteq_{nm}G_{\tilde{C}}$.
First, suppose that $G_{CB}G_{CA}\subseteq_{nm}G_C$. Since $[G_{CB}:G_{\tilde{C}\tilde{B}}],[G_{CA}:G_{\tilde{C}\tilde{A}}]<\omega$, 
we get that $G_{CB}G_{CA}$ is a union of finitely many two-sided translates of $G_{\tilde{C}\tilde{B}}G_{\tilde{C}\tilde{A}}$ by elements of $G_C$. 
So, $G_{\tilde{C}\tilde{B}}G_{\tilde{C}\tilde{A}}$ is non-meager in $G_C$, and, hence, in $G_{\tilde{C}}$.

Now, suppose that $G_{\tilde{C}\tilde{B}}G_{\tilde{C}\tilde{A}}\subseteq_{nm}G_{\tilde{C}}$. Then $G_{CB}G_{CA}\cap G_{\tilde{C}}$ is non-meager in $G_{\tilde{C}}$, and hence, in $G_C$ (because $[G_C:G_{\tilde{C}}]<\omega$). Thus, $G_{CB}G_{CA}$ is non-meager in $G_C$. This proves (1). Now, (2) follows by (1) and transfinite induction. \hfill $\square$\\

The following corollary gives a negative answer to Question \ref{qgeneric} in its full generality, i.e., in the class of all Polish group structures. Recall that Fact \ref{generic} tells us that the answer is positive for small Polish $G$-groups.
\begin{corollary}\label{no generics}
Let $(X,G)$ be a Polish structure, where $X$ is uncountable. If $H$ is a countable group, then $(H(X),G)$ is a small Polish group structure and it has no generic orbit (neither left nor right). 
\end{corollary}
{\em Proof.} Take any $a\in H(X)$ and a finite $A\subseteq H(X)$. We will show that $o(a/A)$ is not a generic orbit.
Take any $d\in X\backslash Acl(\emptyset)$, $h\in H\backslash \{e\}$ and $b\in o(d)$ such that $b\nmind \tilde{A}, \tilde{a}$.
 Then, by Proposition \ref{new ind}, $h_b\nmind_{A} a$. Since $b\nmind \tilde{A}, \tilde{a}$ and $b\notin Acl(\emptyset)$, we see that $b\notin \tilde{a}$. Hence, $\widetilde{a+h_b}=\tilde{a}\cup\{b\}$. 
But $\tilde{a},b\nmdep b$, so, again by Proposition \ref{new ind}, we have that
$a+h_b\nmdep h_b$. 
Hence, $o(a/A)$ is not a generic orbit. \hfill $\square$\\

By the above corollary and Fact \ref{generic}, we get in particular that there is no Polish topology on $H(X)$ such that $(H(X),G)$ is a  $G$-group, i.e., such that the action of $G$ on $H(X)$ is continuous.

By \cite{1} and \cite{4}, every $nm$-stable compact $G$-group is nilpotent-by-finite. When the assumption of compactness is dropped, the corresponding questions concern searching for a subgroup of countable index having some nice algebraic properties. The algebraic structure of $nm$-stable Polish $G$-groups remains unexplored. The following corollary shows that in general, not much can be said about the algebraic structure of small $nm$-stable Polish group structures.
\begin{corollary}
Let $(X,G)$ be an uncountable, small, $nm$-stable Polish structure, and $H$ a non-solvable, countable group. Then $(H(X),G)$ is a small, $nm$-stable Polish group structure, which is not solvable-by-countable.
\end{corollary}
{\em Proof.} By Proposition \ref{new ind}(2), $(H(X),G)$ is $nm$-stable. Now, take a subgroup $A$ of countable index in $H(X)$. Then, there is some $x\in X$, such that $\pi_x[A]=H_x$, where $\pi_x:H(X)\to H_x$ is the projection on the $x$-th coordinate. Since $H_x$ is not solvable, we get that $A$ is not solvable. Thus, $H(X)$ is not solvable-by-countable. \hfill $\square$\\




Now, we will give a variant of the above construction.
Suppose $R$ is a countable commutative ring, and $(X,G)$ is a small Polish structure. Let $R(X)=R[(y_x)_{x\in X}]$ be the ring of polynomials in variables $(y_x)_{x\in X}$ with coefficients in $R$. Then $G$ acts on $R(X)$ by $gw(y_{x_1},\dots,y_{x_n})=w(y_{gx_1},\dots,y_{gx_n})$. If $R$ is a countable field, we can additionally consider $R(X)_0=R((y_x)_{x\in X})$, the field of rational functions in variables $(y_x)_{x\in X}$ with coefficients in $R$. Then, $G$ acts on $R(X)$ by $gf(y_{x_1},\dots,y_{x_n})=f(y_{gx_1},\dots,y_{gx_n})$. As for $H(X)$, one can check that $(R(X),G)$, $(R_0(X),G)$ are small Polish structures. Moreover, if we define $\tilde{w}$ as the set of all $x\in X$ such that $y_x$ occurs in the reduced form of $w$, then we get the same description of $nm$-independence for $(R(X),G)$,$(R_0(X),G)$ as was done for $(H(X),G)$ in Proposition \ref{new ind}. Thus, we get that these structures (which we could call Polish ring structures and Polish field structures) have no generics (in the sense of the additive group), and hence, there is no Polish topology on $(R(X),G)$ or on $(R_0(X),G)$ such that the action of $G$ is continuous.

\section{A non-zero-dimensional small Polish $G$-group}\label{dim}
In this section, we construct a first example of a small non-zero-dimensional Polish $G$-group.

We define a structure of a group on the complete Erd\"{o}s space as in \cite[Proposition 4.3]{3}, but we choose $p=1$ instead of $p=2$ in order to avoid computational complications.
Namely, we let $C\subseteq \mathbb{R}$ be the ternary Cantor set, and $X=C^{\omega}\cap l_1$. By Theorems \ref{criterion} and \ref{erdos}, $X$ (considered with the topology induced from $l_1$) is homeomorphic to the complete Erd\"{o}s space. Consider the standard bijection $\phi :2^{\omega}\to C$ and the product map $\psi:=\phi^{\omega}:(2^{\omega})^{\omega}\to C^{\omega}$. Then it easy to see that $H:=\psi^{-1}[X]$ is a subgroup of $(2^{\omega})^{\omega}$ (we will identify the latter group with $2^{\omega\times\omega}$ in the natural way), and becomes a Polish group with the topology induced from $X$ by $\psi$ (and is homeomorphic to the complete Erd\"{o}s space). This topology is generated by the norm $||z||:=||\psi(z)||_1$, $z\in H$. We also put $||z||=\infty$ if $z\in 2^{\omega\times\omega}\backslash H$. For a subset $A$ of $\omega\times\omega$ we define $||A||:=||\chi _A||$, where $\chi _A$ is the characteristic function of $A$.

Now, we will define an action of a Polish group $G$ on $H$. Let $G_1$ be the group of all permutations of $\omega\times\omega$. 
For any $g\in G_1$ we define the support of $g$ to be $supp(g)=\{a\in \omega\times\omega :g(a)\neq a\}$.  
We put: $$G=\{g\in G_1:||supp(g)||<\infty\}<G_1.$$
It is clear that for any $g\in G$ and $h\in H$, the composition $h\circ g:\omega\times\omega\to 2$ is an element of $H$ (since $||h\circ g||\leq ||h||+||supp(g)||$). Hence, we can define an action of $G$ on $H$ by $gh=h\circ g^{-1}$. Then, $G$ acts on $H$ as automorphisms (both algebraic and topological). Notice, however, that if we consider $G$ with the product topology, then this action is not continuous. Hence, we need another topology on $G$.

We define a metric $d$ on $G$: $$d(f,g)=||supp(f^{-1}g)||.$$ We will consider $G$ with the topology generated by $d$.

\begin{proposition}\label{polish group}
$G$ is a Polish group.
\end{proposition}
{\em Proof.}
It is easy to check that $d$ is a complete metric on $G$. Also, the set of elements of $G$ with finite support is a countable, dense subset of $G$. Now, we will check that the composition $\circ :G\times G\to G$ is continuous.
For any $(f,g),(f_1,g_1)\in G\times G$ we have $d(fg,f_1g_1)=||supp((fg)^{-1}f_1g_1)||=||supp(g^{-1}(f^{-1}f_1g_1g^{-1})g)||$. 
Clearly the composition is continuous at $(e,e)\in G\times G$ (since $supp(fg)\subseteq supp(f)\cup supp(g)$), so it is enough to check that conjugating by $g$ is continuous at $e\in G$. We will check that for every $f\in G$ conjugating by $f^{-1}$ is continuous at $e\in G$, which is of course sufficient.
Notice that $supp(fhf^{-1})=f[supp(h)]$. 
For any $\epsilon >0$ there is $n<\omega$ such that $||supp(f)\backslash n\times\omega||<\epsilon$, 
and since $supp(f)\cap n\times\omega$ is finite,
 we can choose $m<\omega$ such that $f[H\backslash m\times\omega]\subseteq H\backslash n\times\omega$. 
So, for $h$ so close to $e$ that $supp(h)\subseteq H\backslash m\times\omega$, 
we have that 
$f[supp(h)]\subseteq supp(h)\cup (f[supp(h)]\cap supp(f))\subseteq supp(h)\cup(supp(f)\backslash n\times\omega)$. 
This shows that the conjugation by $f$ is continuous, and hence, so is the group composition. 
Similarly one checks that the group inversion on $G$ is continuous.  \hfill $\square$\\

The next proposition shows that we have constructed (the first known) example of a small, non-zero-dimensional Polish $G$-group.

\begin{proposition} \label{main}
$(H,G)$ is a small, Polish $G$-group.
\end{proposition}
{\em Proof.}
To check that the action of $G$ on $H$ is continuous at every $(g,h)\in G\times H$, consider any $(g_1,h_1)\in G\times H$. Then, the functions $gh=h\circ g^{-1}$ and $g_1h_1=h_1\circ g_1^{-1}$ agree on the set $\{a\in \omega\times\omega:g^{-1}(a)=g_1^{-1}(a)\}\cap \{a\in \omega\times\omega:h(g^{-1}(a))=h_1(g^{-1}(a))\}$. The complement of this set in $H$ is the union of $supp(g_1g^{-1})$ and $g[\{a\in \omega\times\omega :h(a)\neq h_1(a)\}]$. For $(g_1,h_1)$ sufficiently close to $(g,h)$ these sets are arbitrary small in the sense of $||\cdot||$ (by a similar argument to the one in the proof of Proposition \ref{polish group}).
So, the action of $G$ on $H$ is continuous. 

It remains to check that for every finite $A\subseteq H$, there are countably many $G_A$-orbits in $H$.
Fix such an $A$. For any $h\in H$ we put $h_0:=h^{-1}[\{0\}], h_1:=h^{-1}[\{1\}]$. Let $B$ be the Boolean algebra generated by the family of sets: $\{a_0,a_1:a\in A\}$, and let $b_0,b_1,\dots,b_n$ be all its atoms. For exactly one $i\leq n$ we have $||b_i||=\infty$, and we assume that this is the case for $i=0$. We will show that the $G_A$-orbit of an element $x$ of $H$ depends only on the cardinalities of the sets $x_0\cap b_0, x_0\cap b_1,\dots, x_0\cap b_n$ and $x_1\cap b_0, x_1\cap b_1,\dots, x_1\cap b_n$. Suppose that for two elements $x,y\in H$ these cardinalities are the same. For all $0\leq i\leq n$ let $g_i$ be a permutation of $b_i$ such that $g_i[b_i\cap x_0]=b_i\cap y_0$ and $g_i[b_i\cap x_1]=b_i\cap y_1$. It is easy to see that $g_0$ can be chosen such that $||supp(g_0)||<\infty$. Then $\bigcup_{i\leq n}g_i$ is an element of $G_A$, and $gx=y$. This completes the proof. \hfill $\square$
\begin{proposition}
$(H,G)$ is not $nm$-stable.
\end{proposition}
{\em Proof.} For $c\in 2^{\omega \times \omega}$, we will write $c_{ij}$ instead of $c(i,j)$. Consider $o=o(a/\emptyset)$, where $a_{ij}=1$ if $j=0$ and $a_{ij}=0$ if $j>0$.
For any $n<\omega$, let $b_n\in H$ be given by $(b_n)_{ij}=1$ if $j=n+1$, and $(b_n)_{ij}=0$ if $j\neq n+1$.
Then, by the proof of Proposition \ref{main}, for every $n<\omega$, we have $$o(a/b_{<n})=\{x\in H: |(i,j)\in\omega \times \omega :x_{ij}=1|=\omega \wedge \forall j
\in\{1,2,\dots n\} x_{ij}=0 \}.$$ So, $o(a/b_{<n})$ is a $G_{\delta}$ subset of $H$. Morever, for every $n<\omega$, $o(a/b_{<n+1})$ is nowhere dense in $o(a/b_{<n})$. Thus, by Fact \ref{gdelta}, $\NM(o)=\infty$. So, $(H,G)$ is not $nm$-stable.  \hfill $\square$\\

Similarly, one can check that $\NM$-rank of every uncountable 1-orbit in $(H,G)$ is equal to $\infty$.
\begin{question}
Is there an $nm$-stable, non-zero-dimensional small Polish $G$-group?
\end{question}

Notice that since the product $H\times H$ is homeomorphic to $H$, we can not obtain examples of higher dimension just by taking finite cartesian powers of $H$. 

\begin{question}
Is there an $nm$-stable, small Polish $G$-group of dimension greater than one?
\end{question}

\noindent
{\bf Address:}\\
Instytut Matematyczny, Uniwersytet Wroc\l awski,\\
pl. Grunwaldzki 2/4, 50-384 Wroc\l aw, Poland.\\[3mm]
{\bf E-mail address:}\\
dobrowol@math.uni.wroc.pl \\

\end{document}